\documentclass{article}

\usepackage{amsmath}
\usepackage{amssymb}
\usepackage{mathpartir}
\usepackage{url}

\newcommand{\Prop}{\mathbf{Prop}}
\newcommand{\False}{\mathrm{False}}
\newcommand{\Prf}{\mathbf{Prf}}
\newcommand{\Type}{\mathbf{Type}}

\newcommand{\Nat}{\mathbb{N}}
\newcommand{\suc}{\mathrm{suc}}
\newcommand{\Int}{\mathbb{Z}}
\newcommand{\Rat}{\mathbb{Q}}
\newcommand{\Real}{\mathbb{R}}
\newcommand{\Compl}{\mathbb{C}}
\newcommand{\Bool}{\mathrm{Bool}}
\newcommand{\true}{\mathrm{true}}
\newcommand{\false}{\mathrm{false}}
\newcommand{\class}{\mathrm{class}}
\newcommand{\refl}{\mathrm{refl}}
\renewcommand{\implies}{\Rightarrow}

\author{Thorsten Altenkirch
}
\title{Should Type Theory replace Set Theory as the Foundation of Mathematics ?}

\begin{document}

\maketitle

\section{Introduction}
\label{sec:introduction}

Set theory is usually traced back to Cantor who used sets in an informal way giving rise to what is called \emph{naive set theory}. Nowadays we usually refer to axiomatic set theory which was formulated by Zermelo and Fraenkel and which is referred to as Zermelo-Fraenkel Set Theory or short ZFC. When saying Set Theory
\footnote{I am capitalising both \textbf{Set Theory} and \textbf{Type Theory} when referring to specific theories with this name not in the sense of a set theory or a type theory.}
we mean ZFC.

Type Theory was introduced by Per Martin-L\"of \cite{martin1975intuitionistic} and there are several incarnations. The early \emph{Extensional Type Theory} (ETT) \cite{martin1984intuitionistic} gave way to \emph{Intensional Type Theory} (ITT) \cite{nordstrom1990programming} but recently, heavily influenced by Voevodsky and concepts from Homotopy Theory, Homotopy Type Theory (HoTT) \cite{hottbook} was developed. 
\footnote{Since not everybody wants to read the whole book, introductory presentations of HoTT include 
\cite{altenkirch2019naive,ahrens2019univalent,grayson2018introduction}.}
When saying Type Theory we mean HoTT.

In Set Theory all mathematical objects are viewed as sets and we write $a \in A$ to mean that $a$ is an element of the set $A$. Since in modern Set Theory we don't use \emph{urelements}, $a$ is a set again, a process that only stops at the empty set $\{\}$. The same element can occur in different sets and two sets are equal iff they have the same elements (axiom of extensionality). 

Set Theory
 is usually based on classical logic, formally it is presented in the framework of first order predicate logic. Set Theory can be viewed as an alternative to higher order logic where the quantification over predicates is replaced by quantifying over sets, corresponding to predicates. There are alternative set theories which use intuitionistic predicate logic as the framework but they all have in common that we first fix the logic as a framework and then formulate the set theory within this logic. 

Type Theory is based on the idea that all mathematical objects belong to a type and can only be understood as elements of a given type. This corresponds to mathematical practice where we conceive of a statement quantifying over all natural numbers as only talking about those while in Set Theory this is represented by quantifying over all sets and singling out those which are elements of a particular one, the set of natural numbers. In Type Theory an element can only be an element of one type
\footnote{I am here presenting a particular pure form of Type Theory. They are extensions which allow subtyping, which is usually an attempt to recover features of Set Theory. I prefer to view these sort of developments as notational conveniences as opposed to foundational features.}
hence we cannot talk about elements in isolation which means we cannot reason about representations. This discipline enables a high degree of extensionality which is realised by Voevodsky's univalence axiom entailing that isomorphic structures are equal.

Types have elements which we write as $a : A$ meaning $a$ is an element of type $A$. While superficially similar to the set-theoretic $a \in A$ it is fundamentally different. Saying $a : A$ is part of the static structure of a mathematical text
\footnote{The static structure of a text corresponds to the static semantics of a program in the sense of computer science. It covers aspects like syntax but also correct use of variables and typing.}
, only texts which use types consistently make sense. In particular we cannot ask whether $a : A$ because this is an aspect of our text not the subject we are talking about. We say that $a : A$ is a judgement not a proposition.

The formal presentation of Type Theory also requires another judgement which is called definitional or judgemental equality. E.g. when we define $n : \Nat$ as $n :\equiv 3$ we also know that $(1,2,3) : \Nat^n$ but to see this we need to exploit the fact that $n \equiv 3$ statically. Again $a \equiv b$ is a judgement which is a static property of a mathematical text  not a predicate. Judgemental equality also reflects the computational character of Type Theory, e.g. we know that all (closed) natural number expressions can be reduced to a numeral, e.g. $3 + 4 \equiv 7$. Type Theory also introduces \emph{propositional equality} but it makes sense to distinguish the static judgemental equality from the usual propositional equality. Indeed being equal by definition is a common phrase in mathematics.


Mathematicians often consider ZFC as the only foundation of mathematics, and frequently don't actually want to think much about foundations. We argue here that modern Type Theory, i.e. HoTT, is a preferable and should be considered as an alternative. 

\subsection{Related work}
\label{sec:related-work}

In \cite{maddy} the author develops a set of criteria for possible foundations and compares set theory, category theory and homotopy Type Theory using these criteria. First of all it seems rather misleading to include category theory as a possible foundation because while it provides \emph{essential guidance} (one of the criteria) it requires another foundation to explain what basic constructions are intuitively permissable.
\footnote{Indeed in loc.cit, Maddy complains: \emph{I don’t know what Voevodsky finds lacking in category-theoretic
foundations -- perhaps that it fails to provide a Generous Arena?}}
 Indeed categorical guidance can be misleading since it would suggest that symmetric cartesian categories would be a good concept but you need to go beyond category theory to see that only preorders satisfy this condition. However, category theory, fits very well with HoTT, indeed the universal properties have unique solutions in HoTT, while they are only unique up to isomorphisms. Maddy in loc.cit., footnote 28,  admits that she is \emph{not sure what these thinkers take to be wrong with ZFC}, where the thinkers seem to refer to Voevodsky and others. It is the purpose of this paper to explain what we think is wrong with ZFC and can be addressed by using HoTT. Maddy also argues that she thinks that HoTT would need to adopt the axiom of choice to be able to encode the classical theory of sets 
\footnote{In loc.cit, p.29: \emph{But to get even to ETCS, we have to add the axiom of choice,} where ETCS stands for \emph{The Elementary Theory of the Category of Sets}}
--- which is in my experience
a common reaction of philosophers who are only too willing to give in to intimidation by classical mathematicians and let them define the yardsticks of mathematics. Her conclusion is that there is no problem with set theory and that the only advantage HoTT would have to offer is proof checking. While it is doubtful that this even correct (there are proof systems based on classical mathematics), it is also misleading. What we try to argue in this paper is that the fundamental point of HoTT is conceptual, the fact that it lends itself to implementation is a welcome additional benefit.

Homotopy Type Theory has attracted some attention in the philosophy community: \cite{awodey2014structuralism} explains why HoTT is important for structural mathematics; \cite{ladyman,ladyman-2} give a philosophical explanation of the identity types in HoTT and \cite{corfield2020modal} how an extension of HoTT namely modal HoTT can be used in natural sciences especially in physics. In the present paper we attempt to argue speccifically what advantages HoTT has over Set Theory. To keep the presentation self-contained we explain the basic ideas of HoTT but in an informal way avoiding the complex formal structure of type theory. A novelty of our paper is that in section \ref{sec:struct-constr} we argue that constructivity and structuralism are just two sides of the same coin.



\section{The role of logic}
\label{sec:role-logic}

Type Theory does not rely on predicate logic as a framework but introduces logic using the propositions as types explanation --- in this sense logic is an emergent aspect of Type Theory not a prerequisite. In constructive thinking we acknowledge that mathematical constructions are taking place in our minds and that we communicate them using shared intuitions. From this point of view it is not the notion of truth which is fundamental but the notion of evidence or proof. A proposition is given by saying what we need as evidence to accept that the proposition holds. In the context of Type Theory we can do this by assigning to every proposition the type of evidence for this proposition. This is also called the Curry-Howard equivalence but this is based on the assumption that we already have a logic and a Type Theory and then observe a formal relationship between the two, while I \emph{define} logic by the translation.
\footnote{The propositions as types explanation is also related to the Brouwer-Heyting-Kolmogorov (BHK) semantics which introduces untyped realisers of propositions.}

The translation of propositional connectives into type theoretic connectives is rather straightforward, for example we translate implication $P \implies Q$ as the function type $P \to Q$ between the associated types of evidence which I denote by the same letters here
\footnote{To be more precise we can introduce an operation $\mathrm{Prf}$ which assign to every proposition $P$ te type of its proofs $\Prf(P)$ and then we stipulate that $\Prf(P \implies Q) \equiv \Prf(P) \to \Prf(Q)$.}
. The idea is that evidence for $P$ implies $Q$ is a function which maps evidence for $P$ into evidence of $Q$ --- function types are a primitive concept of Type Theory
\footnote{We explain a function as \emph{black box} where we can input elements of the domain and which outputs elements of the codomain. While we have no access to the mechanism (hence the box is black) the only way to actually construct a function is via effective means. However, the theory works as well if we believe in some sort of magic, like that everything is decidable (the principle of the excluded middle). Formally, the properties of functions are represented by the laws of $\Lambda$-calculus and the principle of functional extensionality.}.
Similar we can translate $\forall x:A.P(x)$ by the dependent function type
\footnote{The dependent function type is a generalisation of the usual function type where the codomain can vary over the domain. This corresponds to the cartesian product of a family of sets in Set Theory, hence the use of $\Pi$}
denoted as $\Pi x:A.P(x)$, the idea is that evidence for a universally quantified statement is a (dependent) function that assigns to every element $a : A$ evidence for $P(a)$. Negation $\neg P$ we interpret as saying that $P$ is impossible means $P$ implies falsum ($P \implies \False$) where falsum ($\False$) is interpreted as the empty type.

The translation of disjunction and existential quantification depends on what exactly we mean by a proposition. Here there has been a historic change mainly since the introduction of HoTT. In ITT we identified types and propositions which leads to some strange artefacts. For example we would denote a subset $\{ x : A \mid P(x) \}$ as the type of dependent pairs $(a,p)$ where $a : A$ and $p : P(a)$ (this is written as $\Sigma x:A.P(x)$)
\footnote{The justification for this notation is similar as before, in Set Theory the type of dependent pairs corresponds to an infinite sum hence $\Sigma$.}
. Now the embedding from the subset into $A$ is given by the first projection $\pi_1  : \{ x : A \mid P(x) \} \to A$
\footnote{Here $\pi_1 : \Sigma x:A.B(x) \to A$ is defined as $\pi_1(a,b) \equiv a$ and $\pi_2 : \Pi p : \Sigma x:A.B(x). B(\pi_1(p))$ with $\pi_2(a,b)\equiv b$.}.
However, it turns out that this is not in general an injection since we can have different proofs $p$ for the same element $a$. Also the naive translation of existential statements uses again dependent pairs, i.e. evidence for $\exists x:A.P(x)$ is given by a pair $(a,p)$ with $a : A$ and $p : P(a)$, we write this type as $\Sigma x:A.P(x)$ as before. However, this interpretation of existence gives rise to an apparent proof of the axiom of choice: if we know that for all $x : A$ there exists a $y : B$ such that a certain relation $R(x,y)$ holds then we can extract the witness and its justification from the proof. More precisely we can define
\begin{align*}
  &\mathrm{ac} : (\forall x : A.\exists y:B . R(x,y)) \to \exists f : A \to B . \forall x:A.R(x,f(x))\\
  &\mathrm{ac}(f) = (\lambda a. \pi_1(f(a)), \lambda a.\pi_2(f(x)))
\end{align*}
However, this \emph{type theoretic axiom of choice} hardly corresponds to the axiom of choice as it is used in Set Theory. Indeed, it is not an axiom but just a derivable fact.

In HoTT we say a proposition is a type with at most one element. The idea here is that a proposition should not contain any additional information, its only purpose is to assert something. As a consequence of univalence we obtain that two propositions are equal if they are logically equivalent --- this is often called propositional extensionality. The embedding from a subset to a set ($\pi_1  : \{ x : A \mid P(x) \} \to A $) is now an injection because there is at most one proof that the proposition $P(a)$ holds. To translate disjunction and existence we use propositional truncation, that is to any type $A$ we introduce a proposition $||A||$ which corresponds to the proposition that $A$ is inhabited.
\footnote{Formally, $||A||$ can be defined as a higher inductive type with two constructors $|-| : A \to ||A||$ and $\mathrm{is-prop} : \forall x y : ||A||. x= y$.}
That is if there is an element $a:A$ we obtain an element of $||A||$ but we cannot distinguish them anymore.  Now we can say that a proof of existence, i.e. evidence for $\exists x:A.P(x)$ is the propositional truncation of the naive translation, i.e. $||\Sigma x:A.P(x)||$. Given this interpretation the translation of
\[ (\forall x : A.\exists y:B . R(x,y)) \to \exists f : A \to B . \forall x:A.R(x,f(x))\]
is not provable since we extract an element from a propositional truncation which is clearly wrong because here all elements have been identified. And indeed this formulation of the axiom of choice has the power of the classical axiom, in particular we can derive the excluded middle using Diaconescu's construction.

By the principle of the excluded middle
\footnote{The term \emph{excluded middle} doesn't actually reflect its meaning. \emph{Excluded} clearly is a negative statement, hence we should render it is saying \emph{It cannot be that a proposition is neither true nor false}, which can be translated is $\neg\neg(P \vee \neg P)$, which \emph{is} constructively provable. A better term for the universal assertion of $P\vee\neg P$, would be \emph{universal decidability} or \emph{magicae nigrae} in latin.} 
we mean
\[ \forall P : \Prop . P \vee \neg P \]
Here I write $\Prop$ for all \emph{propositional types}, i.e. types with at most one element and $P \vee Q$ is the translation of disjunction using propositional truncation, ie. $|| P + Q||$
\footnote{Here $A+B$ is the coproduct or binary sum in Type Theory corresponding to disjoint union in Set Theory. It can be defined as $A+B = \Sigma b : \Bool.P(B)$ with $P(\true) \equiv A$ and $P(\false)\equiv B$.}
. Actually in this particular example we don't need to truncate because it can never happen that both $P$ and $\neg P$ hold. Clearly there is no function which assigns to every proposition either a proof of it or a proof of its negation. We don't even need to appeal to the undecidability of our logic but observe that we cannot look \emph{into} a proposition due to propositional extensionality. Hence excluded middle seems to suggest that we can decide any proposition without even looking at it!

However, we can justify excluded middle by a different translation of the connectives and a different explanation what is a proposition. That is we can interpret $P \vee_\class Q$ as the statement that not both are false $\neg(\neg P \wedge \neg Q)$ and similar we can explain existence $\exists_\class x:A.P(x)$ by saying it is not the case that $P$ is always false, i.e. $\neg (\forall x:A.\neg P(x))$. Note that we don't need any truncation because negated statements are always propositions. Moreover the translations are valid in classical logic anyway. However, to obtain the correct behaviour of the connectives, in particular to justify reasoning by cases and the elimination rule of the existential, we need to limit propositions to the negative  ones, i.e. to the ones for which proof by contradiction holds, that is we define $\Prop_{\mathrm class} = \{ P : \Prop \mid \neg \neg P \to P\}$. So for example if we know $P\vee Q$ and $P \to R$ and $Q \to R$ we can conclude $\neg \neg R$ and if $R$ is negative we also know $R$. Similar reasoning justifies the rule of existential elimination. Moreover the remaining connectives preserve the property of being negative. What I present here is just the usual negative translation.

We may think that there is no difference between propositions as types with at most one element and negative propositions because in neither case can we extract information. However, while the axiom of choice is unprovable in both cases, the axiom of unique choice is provable for the type theoretic interpretation of propositions, that is
\begin{align*}
  &\mathrm{auc} : (\forall x : A.\exists! y:B . R(x,y)) \to \exists f : A \to B . \forall x:A.R(x,f(x))\\
  &\mathrm{auc}(f) = (\lambda a. \pi_1(f(a)), \lambda a.\pi_2(f(x)))
\end{align*}
where $\exists! x:A.P(x)$ means \emph{unique existence} which can be defined as
\[\exists x:A.P(x)\wedge \forall y:A.P(y) \implies x=y.\]
This can be proven in Type Theory because the statement of unique existence is already a proposition anyway and hence we can extract the witness. However, using the negative translation we cannot justify unique choice. 

\section{Representation}
\label{sec:representation}

To illustrate the difference between Set Theory and Type Theory let's look at the representation of numbers, starting with the natural numbers. In Set Theory natural numbers are encoded using only sets, the most common encoding is due to von Neumann:
\begin{align*}
  0 & = \{\} \\
  1 & = \{ 0 \}\\
  2 & = \{ 0 , 1 \}\\
      \vdots & \quad \vdots \\
  n & = \{ 0 , 1, \dots , n-1 \}
\end{align*}
However, there are many alternative ways, e.g. the following is due to Zermelo:
\begin{align*}
  0 & = \{\} \\
  1 & = \{ 0 \}\\
  2 & = \{ 1 \}\\
      \vdots & \quad \vdots \\
  n & = \{ n-1 \}
\end{align*}
So for example the number 3 is represented as $\{\{\},\{\{\}\},\{\{\},\{\{\}\}\}$ using von Neumann's encoding and as $\{\{\{\}\}\}$ in Zermelo's. The axiom of infinity uses one particular representation to introduce the existence of the infinite set of natural numbers, e.g. von Neumann's, and using the axiom of replacement we can deduce that the alternative representation, e.g. Zermelo's also forms a set.

We can formulate properties that distinguish the two representations, for example in von Neumann's encoding we have that $n \subseteq n+1$ which fails in Zermelo's, e.g. $1 \not\subseteq 2$ if we use Zermelo's encoding. This statement is not really a statement about numbers but about their encoding, we say it is not structural. Mathematicians would normally avoid non-structural properties, because they entail that results are may not be transferable between different representations of the same concept. However, frequently non-structural properties are exploited to prove structural properties and then it is not clear whether the result is transferable.

In Type Theory natural numbers are an example of an inductive type, defined by the constructors. Following Peano we say that natural numbers are generated from $0 : \Nat$ and $\suc : \Nat \to \Nat$ (the successor, $\suc(n)$ means $n+1$). To define functions out of the natural numbers we use an elimination principle that allows us to perform dependent recursion, that is given a family of types indexed by the natural numbers $M : \Nat \to \Type$ and an element $z : M(0)$ and $s : \Pi n:\Nat . M(n) \to M(\suc(n))$ we can define a function $f : \Pi n:\Nat.M(n)$ with the definitional equalities $f(0)\equiv m(0)$ and $f(\suc(m)) \equiv s(n,f(n))$. The elimination principle allows us to define function like addition but also to derive properties, e.g. that addition forms a commutative monoid exploiting the propositions as types explanation.

The definition of the natural numbers precisely identifies the structural properties of the type: it is equivalent to stating that $\Nat$ is an initial algebra, a natural number object, in the language of category theory. And because we cannot talk about elements in isolation it is not possible to even state non-structural properties of the natural numbers. Indeed, we cannot distinguish different representations, for example using binary numbers instead. 

Both in Set Theory and in Type Theory we can go on and construct more elaborate sets/types of numbers, e.g. $\Int,\Rat,\Real,\Compl$. In Set Theory we can talk about subsets, e.g. we may want to stipulate that $\Nat \subseteq \Int$ but it turns out that this doesn't hold for the most natural choice of encodings.
\footnote{E.g. we may represent integers as equivalence classes of pairs of natural numbers with the same difference, i.e. $(a,b) = (c,d)$ if $a+d = c+b$.}.
Yes, we can start with the biggest number class we want to represent and then define the other ones as subsets via comprehension.

In Type Theory the notion of a subset doesn't make sense in general for types. We can observe that we have a canonical embedding function between the different number types, e.g. $\Nat \to \Int$ and we can agree notationally that we can omit this function when it is clear from the context that we use a natural number when an integer was required.

While we don't have a notion of subset on types, we can talk about subsets of a given type $A$ which are just propositionally valued functions $P : A \to \Prop$. Given any two of such subsets we can define $P \subseteq Q$ as $\forall x:A.P(x) \to Q(x)$ in Type Theory. This means we can actually play the same trick as in Set Theory and define our number classes as subsets of the largest number class we want to consider and we have indeed the subset relations we may expect.

Hence Type Theory allows us to do basically the same things as Set Theory as far as numbers are concerned (modulo the question of constructivity) but in a more disciplined fashion limiting the statements we can express and prove to purely structural ones.

One criteria for foundation might be minimality: we need as few basic constructs as possible. One may think at the first glance that Type Theory is somehow less minimalistic, while in Set Theory we only need $\{\dots\}$ to represent all mathematical concepts there seems to be quite a menagerie of type constructors: $\Pi$-types, $\Sigma$-types, equality types and so on. However, to be fair we should compare the axioms of Set Theory with the basic constructions of Type Theory. And indeed while in applied Type Theory we may introduce a lot of apparently unrelated types we can often reduce them to a basic collection. So for example most inductive constructions can be reduced to the type of well-founded trees or $W$-types. However, it s fair to admit that the program to reduce very advanced type theoretic constructions, e.g. Higher Inductive Types (HITs) to basic combinators is subject of current research. This can be compared to reducing advanced set theoretic principle to basic axioms. 

\section{Univalence and equality}
\label{sec:univalence-equality}

Type Theory doesn't allow us to make statements about representations because we cannot talk about elements in isolation. This means that we cannot observe intensional properties of our constructions. This already applies to Intensional Type Theory, so for example we cannot observe any difference between two functions which are pointwise equal. However, on the other hand in ITT we cannot prove them to be equal either. This is a consequence of the way the equality type is treated in ITT. That is for any element $a : A$ there is a canonical element of the equality type $\refl(a) : a = a$. Now consider $f,g : \Nat \to \Nat$ with
\begin{align*}
  f(x) & = x + 0 \\
  g(x) & = 0 + x
\end{align*}
Using induction it is easy to show that $f(x) = g(x)$ for all $x:\Nat$. However, the only prove of $f = g$ would need to be canonical (since we have no other assumptions) i.e. $\refl(h)$ but this would require that $h \equiv f \equiv g$ and $f$ and $g$ are clearly not definitionally equal.

Hence in ITT while we cannot distinguish extensionally equal functions we do not identify them either. This seems to be a rather inconvenient incompleteness of ITT, which is overcome by HoTT. In HoTT a type doesn't just have elements but also it comes with an equality type that tells what is the type of evidence that two objects are equal. In the case of function type this evidence is equivalent to functional extensionality. Another example is propositional extensionality: two propositions are equal if and only if they are logically equivalent. By incorporating these principles Type Theory is as extensional as Set Theory but it turns out it is much better.

This is a consequence of the fact that we cannot observe the representations of objects. What does this mean for types: what can we observe about a type? Now we can observe the cardinality of a type, e.g. whether it has two elements or whether it is countably infinite but  nothing else. Hence extensionally two types with the same number of elements should be considered equal. This may appear counterintuitive on the first glance because it entails for example that $\Int = \Nat$. But then can you write a property of a type that distinguished $\Int$ from $\Nat$? You may say that addition behaves differently for the two but you just have the type and no operation on them. Indeed, if we talk about the algebraic structure, e.g. we consider them as commutative monoids, then they are not equivalent because we can distinguish them.

It is a consequence of the univalence principle that two isomorphic types are equal
\footnote{The actual statement of univalence is slightly different and uses the notion of an \emph{equivalence} which refines isomorphism in a setting where equality is not always propositional. A function $f : A \to B$ is an equivalence if for any $b : B$ there is a unique $a :A$ \textbf{and} $p : f(a) = b$. Since equivalence is reflexive there is an obvious function from equality of two types to their equivalence. Univalence states that this function is an equivalence. }
where by isomorphic I mean that there a functions in both directions which when composed are equal to the identity. So basically two types are equal if there is a one-to-one correspondence between their elements. This also applies to algebraic structures like the aforementioned commutative monoids. While the monoids associated to $\Int$ and $\Nat$ are not isomorphic we cannot distinguish different representations of the same algebraic structure, e.g. using unary or binary numbers to represent $\Nat$. And indeed it is a consequence of univalence that isomorphic algebraic structures are equal.

This reflects mathematical practice to view isomorphic structures as equal. However, this is certainly not supported by Set Theory which can distinguish isomorphic structures. Yes, indeed all structural properties are preserved but what exactly are those. In HoTT all properties are structural, hence the problem disappears.

There is a price to pay for this elegance. We usually expect that equality is a proposition but this ceases to hold when we say that evidence for equality of types is isomorphism because there is usually more than one isomorphism between isomorphic types. 
\footnote{One can force equality to be a proposition by truncating it. However, this means that we loose important properties such as the ability to always replace equals by equals.}
The conclusion is that equality of structures in general is not a proposition but a structure itself.

What is equality actually? I said that every type comes with a notion of equality, with an equality type but what are the properties of equality types. If equality were just a proposition then it would be easy: equality is just an equivalence relation (i.e. reflexive, symmetric and transitive) and every function has to preserve it. This view is called the setoid interpretation and it is a useful easy case of the actual story. 
\footnote{Setoids have been used to justify extensional principles, see \cite{altenkirch1999extensional,altenkirch2019setoid}.}
But since we have types and we can iterate the construction of equality types, i.e. we get a whole infinite tower of equalities, the answer is a bit more involved. The first step up from equivalence relations are groupoids, these are categories where every morphism corresponding to an equality type has an inverse. Or more naively, it is the combination of an equivalence relation and a group.
\footnote{A one object groupoid is a group while a groupoid where every homset is propositional is an equivalence relation.}
However, we need to iterate this process any number of times arriving at weak infinity groupoids. They are called weak because all the properties don't hold strictly but only upto higher equality. 
\footnote{\cite{ladyman} provide a philosophical justification of path induction. However, they seem to miss the point that path induction is nothing but a way to express that equality is a higher groupoid, that is the natural generalisation of an equivalence relation to a proof-relevant seeting.}
This is where Homotopy theory comes in helpful because these structures have been investigated and they can be defined as simplicial sets with the Kan property or more exactly Kan fibrations for dependent types. Actually in the context of Type Theory it turns out that cubical sets are technically better behaved.
\footnote{In \cite{cohen2016cubical} cubical sets are used to provide a constructive explanation of the univalence axiom. This semantics gives rise to cubical Type Theory,
  which is now implemented in cubical agda.\cite{vezzosi2021cubical}. Cubical Type Theory can be understood as a refinement of HoTT. }

\section{Structuralism and constructivity}
\label{sec:struct-constr}

We may think that the constructive nature of Type Theory which comes from the propositions as types approach to logic is rather orthogonal to its support of structural mathematics which stems from the fact that elementhood is treated statically. Indeed, we can add classical principles like the unrestricted axiom of choice (and hence excluded middle) to HoTT and still have univalence. So are these two aspects orthogonal, can we just have univalence but don't embrace constructive Mathematics.

While this is possible in principle, I would argue that there is indeed a connection between these two aspects. Univalence forces us to revisit the notion of equality and accept that equality isn't always a proposition. It is this step which at the same time renders many applications of choice unnecessary. Choice in conventional mathematics becomes necessary because we only talk about propositions, and hence choices need to be hidden in propositions. When we later need to refer to the choices we have just hidden, we need to refer to the axiom of choice. In Type Theory we can indeed exploit the fact that our constructions may be types and not propositions and the choices are explicit.

A simple example for this phenomena is the way categories are treated in HoTT. Categories are a generalisation of preorders but indeed we prefer to deal with partial orders where the equivalence introduced by the order coincides with equality. Can we do the same for categories? Indeed in HoTT we can introduce the notion of a \emph{univalent category}
\footnote{In the terminology of the HoTT book, ordinary categories are called precategories and univalent categories are called categories.}
, where isomorphism and equality of objects are equivalent (i.e. isomorphic). And indeed most categories that arise from semantic consideration -- like the category of sets, or categories of algebraic structure -- are univalent. In contrast the counterpart to univalent categories in a classical setting are \emph{skeletal categories}, which are very rare and their constructions only succeeds by heavy use of the axiom of choice.

\section{Conclusions}
\label{sec:conclusions}

Axiomatic Set Theory was developed in the 1930ies and then it constituted considerable progress in making the foundations of mathematics precise. However, much has happened since and it is time to reconsider the foundations of mathematics, in particular in light of a shift of major applications of mathematics from natural science to computer science and related disciplines. 

From the viewpoint of computer science, it seems almost preposterous to think that all mathematical concepts can be adquately represented as sets, which correspond to a type of trees where we ignore order and multiplicity of subtrees. This encoding seems to be not very different than the insight that all structures in computing can be represented as a sequence of 0s and 1s. Indeed, it seems to me that set theory is something like the machine language of mathematics, but it doesn't share the redeeming feature that in computing we \emph{need} to translate everything in machine language. Mathematics which is a construction of our mind doesn't need a machine language it should rather reflect our basic intuition. 

Set Theory enables us to ask stupid questions, like \emph{Are the integers a subset of the rationals? }or \emph{Is 2 a subset of 3}? Mathematicians tell me that they just don't ask these questions so they don't really need a better language. This reminds me, who worked as a programmer in the 1980ies, of the attitude of some of my colleagues who didn't see the need for high-level languages because they knew how to do the right thing on the machine level.

The basic tenants of Type Theory are this: mathematical concepts are represented as types, where every type comes with a static notion of elementhood and an equality type. Objects are equal if they represent the same mathematical object, in particular we view objects as equal if they satisfy the same properties -- this is just paraphrasing Leibniz
\footnote{Here I disagree with the analysis in \cite{ladyman-2} who look for an internal representation of Leibniz's principle which naturally fails. We need to look atthe philosophical idea of identifying an object with its properties which leads to the notion of extensionality which is realized in HoTT.}.
A consequence of this is that equalities aren't always propositions because they may be me more than one way that two objects can be equal. Logic emerges in Type Theory, as the propositions as types explanation, which replaces the philosophically difficult concepts of truth with the notion of evidence which can be expressed as a type.

As far as constructivity is concerned, I would like to quote Frank Pfenning who recently said in a lecture \emph{I have lost the ability to understand classical logic}. I share the same problem, specifically with the boolean explanation of logic: it seems to me that we cannot \emph{assign} a boolean to every formula in arithmetic. Now I understand that \emph{assign} here seems to refer to a \emph{classical function}. However, and this may be because I am computer scientist, it seems to me a strange abuse of the word assignment if you have one but cannot tell me what is the output for a given input. 

It seems to me that the classical approach can only succeed by an intensive process of brain-washing, replacing naive concepts, like function, which rather elaborate replacements. It seems that many mathematicians think that this process is necessary for them to continue their leading edge research in modern mathematics. That may well be so, but maybe we shouldn't give up the hope that there is the possibility for a more naive form of mathematics which favours explainability over latest fashion research.

\bibliographystyle{alpha}
\bibliography{book,local}

\end{document}